\documentclass{amsart}

\usepackage{amsmath,amsthm,amsfonts,amscd,amssymb}
\usepackage[pdfborder={0 0 0}]{hyperref}
\usepackage{tikz}
\usepackage{float}
\usepackage{amsbsy}

\newcommand{\abs}[1]{\left\lvert #1 \right\rvert}

\newcommand{\ZZ}{\mathbb{Z}}

\newtheorem{thm}{Theorem}[section]
\newtheorem{cor}[thm]{Corollary}
\newtheorem{lem}[thm]{Lemma}

\newtheorem{quest}{Question}

\theoremstyle{definition}

\theoremstyle{remark}
\newtheorem{rem}{Remark}[section]

\numberwithin{equation}{section}

\begin{document}

\title{On well-rounded sublattices of the hexagonal lattice}
\author[L. Fukshansky]{Lenny Fukshansky}
\author[D. Moore]{Daniel Moore} 
\author[R. A. Ohana]{R. Andrew Ohana}
\author[W. Zeldow]{Whitney Zeldow}\thanks{The authors were supported by NSF grant DMS-0755540}

\address{Department of Mathematics, 850 Columbia Avenue, Claremont McKenna College, Claremont, CA 91711}
\email{lenny@cmc.edu}
\address{Department of Mathematics, Loyola Marymount University, 1 LMU Drive, Los Angeles, CA 90045}
\email{daniel.ross.moore@gmail.com}
\address{Department of Mathematics, University of Washington, Box 354350, Seattle, WA 98195-4350}
\email{andrew.ohana@gmail.com}
\address{Department of Mathematics, San Francisco State University, 1600 Holloway Avenue, San Francisco, CA 94132}
\email{whitneyzeldow@hotmail.com}

\begin{abstract} 
We produce an explicit parameterization of well-rounded sublattices of the hexagonal lattice in the plane, splitting them into similarity classes. We use this parameterization to study the number, the greatest minimal norm, and the highest signal-to-noise ratio of well-rounded sublattices of the hexagonal lattice of a fixed index. This investigation parallels earlier work by Bernstein, Sloane, and Wright where similar questions were addressed on the space of {\it all} sublattices of the hexagonal lattice. Our restriction is motivated by the importance of well-rounded lattices for discrete optimization problems. Finally, we also discuss the existence of a natural combinatorial structure on the set of similarity classes of well-rounded sublattices of the hexagonal lattice, induced by the action of a certain matrix monoid.
\end{abstract}

\maketitle
\tableofcontents

\def\A{{\mathcal A}}
\def\AA{{\mathfrak A}}
\def\B{{\mathcal B}}
\def\C{{\mathcal C}}
\def\D{{\mathcal D}}
\def\EE{{\mathfrak E}}
\def\F{{\mathcal F}}
\def\x{{\mathcal H}}
\def\I{{\mathcal I}}
\def\II{{\mathfrak I}}
\def\J{{\mathcal J}}
\def\K{{\mathcal K}}
\def\kk{{\mathfrak K}}
\def\L{{\mathcal L}}
\def\LL{{\mathfrak L}}
\def\M{{\mathcal M}}
\def\mm{{\mathfrak m}}
\def\MM{{\mathfrak M}}
\def\N{{\mathcal N}}
\def\OO{{\mathfrak O}}
\def\PP{{\mathfrak P}}
\def\R{{\mathcal R}}
\def\PNR{{\mathcal P_N(\real)}}
\def\PMNR{{\mathcal P^M_N(\real)}}
\def\PdNR{{\mathcal P^d_N(\real)}}
\def\s{{\mathcal S}}
\def\V{{\mathcal V}}
\def\X{{\mathcal X}}
\def\Y{{\mathcal Y}}
\def\Z{{\mathcal Z}}
\def\H{{\mathcal H}}
\def\cee{{\mathbb C}}
\def\pee{{\mathbb P}}
\def\que{{\mathbb Q}}
\def\QQ{{\mathbb Q}}
\def\real{{\mathbb R}}
\def\RR{{\mathbb R}}
\def\zed{{\mathbb Z}}
\def\ZZ{{\mathbb Z}}
\def\aaa{{\mathbb A}}
\def\ff{{\mathbb F}}
\def\kk{{\mathfrak K}}
\def\qbar{{\overline{\mathbb Q}}}
\def\kbar{{\overline{K}}}
\def\ybar{{\overline{Y}}}
\def\kkbar{{\overline{\mathfrak K}}}
\def\ubar{{\overline{U}}}
\def\eps{{\varepsilon}}
\def\ahat{{\hat \alpha}}
\def\bhat{{\hat \beta}}
\def\gt{{\tilde \gamma}}
\def\h{{\tfrac12}}
\def\be{{\boldsymbol e}}
\def\bei{{\boldsymbol e_i}}
\def\bc{{\boldsymbol c}}
\def\bm{{\boldsymbol m}}
\def\bk{{\boldsymbol k}}
\def\bi{{\boldsymbol i}}
\def\bl{{\boldsymbol l}}
\def\bq{{\boldsymbol q}}
\def\bu{{\boldsymbol u}}
\def\bt{{\boldsymbol t}}
\def\bs{{\boldsymbol s}}
\def\bv{{\boldsymbol v}}
\def\bw{{\boldsymbol w}}
\def\bx{{\boldsymbol x}}
\def\bX{{\boldsymbol X}}
\def\bz{{\boldsymbol z}}
\def\bwy{{\boldsymbol y}}
\def\bY{{\boldsymbol Y}}
\def\bL{{\boldsymbol L}}
\def\ba{{\boldsymbol a}}
\def\bb{{\boldsymbol b}}
\def\bet{{\boldsymbol\eta}}
\def\bxi{{\boldsymbol\xi}}
\def\bo{{\boldsymbol 0}}
\def\bone{{\boldsymbol 1}}
\def\bol{{\boldsymbol 1}_L}
\def\ep{\varepsilon}
\def\p{\boldsymbol\varphi}
\def\q{\boldsymbol\psi}
\def\rank{\operatorname{rank}}
\def\aut{\operatorname{Aut}}
\def\lcm{\operatorname{lcm}}
\def\sgn{\operatorname{sgn}}
\def\spn{\operatorname{span}}
\def\md{\operatorname{mod}}
\def\Norm{\operatorname{Norm}}
\def\dim{\operatorname{dim}}
\def\det{\operatorname{det}}
\def\Vol{\operatorname{Vol}}
\def\rk{\operatorname{rk}}
\def\ord{\operatorname{ord}}
\def\ker{\operatorname{ker}}
\def\div{\operatorname{div}}
\def\Gal{\operatorname{Gal}}
\def\GL{\operatorname{GL}}
\def\SNR{\operatorname{SNR}}
\def\WR{\operatorname{WR}}
\def\scg{\operatorname{\left< \Gamma \right>}}
\def\swrh{\operatorname{Sim_{WR}(\Lambda_h)}}
\def\ch{\operatorname{C_h}}
\def\cht{\operatorname{C_h(\theta)}}
\def\scgt{\operatorname{\left< \Gamma_{\theta} \right>}}
\def\scgmn{\operatorname{\left< \Gamma_{m,n} \right>}}
\def\gat{\operatorname{\Omega_{\theta}}}

\section{Introduction}
\label{intro}

Throughout this paper, we will write $\Lambda_h$ for the hexagonal lattice in $\real^2$:
\begin{equation}
\label{hex_def}
\Lambda_h := \begin{bmatrix} 1&-\frac{1}{2} \\ 0&\frac{\sqrt{3}}{2} \end{bmatrix} \zed^2,
\end{equation}
which is routinely identified with the ring of Eisenstein integers $\zed[\omega]$, where $\omega = e^{\frac{2\pi i}{3}}$. This lattice has many important properties, in particular it solves a variety of classical discrete optimization problems in the plane, such as circle packing and covering problems, kissing number problem, and quantizer problem (see \cite{conway} for a detailed account). Not surprisingly, the properties and structure of $\Lambda_h$ have been extensively studied for their own sake, as well as for the benefit of many applications arising in engineering and digital communications problems. In particular, a detailed analysis of distribution and optimization properties of sublattices of $\Lambda_h$ has been carried out by Bernstein, Sloane, and Wright in \cite{sloane}. It is the goal of this note to continue this investigation, concentrating on the more special class of {\it well-rounded} (from now on abbreviated as WR) sublattices of $\Lambda_h$. 

Given a lattice $\Gamma =A\zed^2 \subset \real^2$, where $A$ is a basis matrix, we define its {\it determinant} to be $\det(\Gamma) = |\det(A)|$, which does not depend on the choice of a basis, and its {\it minimum} (or {\it minimal norm}) to be
$$|\Gamma| = \min \{ \|\bwy\|^2 : \bwy \in \Gamma \setminus \{\bo\} \},$$
where $\|\ \|$ stands for the usual Euclidean norm. Then each $\bx \in \Gamma$ such that $\|\bx\|^2 = |\Gamma|$ is called a {\it minimal vector}, and the set of minimal vectors of $\Gamma$ is denoted by $S(\Gamma)$. A lattice $\Gamma \subset \real^2$ is called WR if there exists $\bx,\bwy \in S(\Gamma)$ which form a basis for $\Gamma$, in which case we call $\bx,\bwy$ a {\it minimal basis}. This minimal basis is not unique, but it is always possible to select a minimal basis $\bx,\bwy$ for a WR lattice $\Gamma$ so that the angle $\theta$ between these two vectors lies in the interval $[\pi/3,\pi/2]$, and any value of the angle in this interval is possible. From now on when we talk about a minimal basis for a WR lattice in the plane, we will always mean such a choice. Then the angle between minimal basis vectors is an invariant of the lattice, and we call it the {\it angle of the lattice} $\Gamma$, denoted $\theta(\Gamma)$; in other words, if $\bx,\bwy$ is any minimal basis for $\Gamma$ and $\gamma$ is the angle between $\bx$ and $\bwy$, then $\gamma = \theta(\Gamma)$ (see \cite{me:hex} for details and proofs of the basic properties of WR lattices in $\real^2$). WR lattices are important in coding theory \cite{esm} and discrete optimization problems \cite{martinet}; they also come up in the context of some number theoretic problems, such as Minkowski's conjecture \cite{mcmullen} and the linear Diophantine problem of Frobenius \cite{me:sinai}. The distribution of WR sublattices of $\zed^2$ has been studied in \cite{me:wr1} and \cite{me:wr}. 

In \cite{sloane}, the authors consider sublattices of $\Lambda_h$ of fixed index $J \geq 2$, counting their number (up to {\it similarity} - to be defined below), and asking which of them have the largest minimum and {\it signal-to-noise ratio}, abbreviated SNR (to be defined in \eqref{SNR_def} below). They provide only a partial answer for these last two questions, proving that both of these quantities are maximized by an {\it ideal sublattice} (i.e., a sublattice coming from an ideal in the ring of Eisenstein integers) whenever there exists one of index $J$. However ideal sublattices, which are a special case of WR sublattices of $\Lambda_h$, do not exist for all possible values of the index, and the authors in \cite{sloane} remark that ``for other values of the index there does not seem to be any general rule to identify which sublattices are best." This motivates a closer investigation of WR sublattices of $\Lambda_h$, which, as we show, exist for more values of the index than ideal sublattices. In this paper we discuss WR sublattices of $\Lambda_h$, giving an explicit description and parameterization for all of them. We then use this parameterization to study their properties with a view toward the discrete optimization questions analogous to those asked in \cite{sloane}.

In order to state our results the notion of similarity of lattices is needed. Two lattices $\Gamma_1,\Gamma_2 \subset \real^2$ are called {\it similar}, denoted $\Gamma_1 \sim \Gamma_2$, if there exists a nonzero real number $\alpha$ and a $2 \times 2$ real orthogonal matrix $A$ such that $\Gamma_2 = \alpha A \Gamma_1$. Similarity is easily seen to be an equivalence relation, and we refer to the equivalence classes under this relation as similarity classes of lattices. WR lattices can only be similar to WR lattices, hence it makes sense to talk about similarity classes of WR lattices. In fact, it is easy to notice that two WR lattices $\Gamma_1,\Gamma_2 \subset \real^2$ are similar if and only if $\theta(\Gamma_1)=\theta(\Gamma_2)$ (see \cite{me:hex} for a proof). Therefore the set of all similarity classes of WR lattices is bijectively parameterized by the set of all possible values of the angle, which is the interval $[\pi/3,\pi/2]$. On the other hand, this parameterization becomes much less trivial if we talk about similarity classes of WR sublattices of $\Lambda_h$. Let us write $\WR(\Lambda_h)$ for the set of all WR sublattices of $\Lambda_h$, and for each $\Gamma \in \WR(\Lambda_h)$ define
\begin{equation}
\label{scg}
\scg := \left\{ \Omega \in \WR(\Lambda_h) : \Omega \sim \Gamma \right\} = \left\{ \Omega \in \WR(\Lambda_h) : \theta(\Omega) = \theta(\Gamma) \right\},
\end{equation}
and let 
$$\swrh := \left\{ \scg : \Gamma \in \WR(\Lambda_h) \right\}.$$
Then it is clear that $\swrh$ is bijectively parameterized by some subset of the interval $[\pi/3,\pi/2]$, and the natural question is what is this subset? Our main result answers this question in detail.

\begin{thm} \label{sim_par} Let 
\begin{equation}
\label{set_ch}
C_h = \{ \theta \in [\pi/3,\pi/2] : \theta = \theta(\Gamma) \text{ for some } \Gamma \in \WR(\Lambda_h) \}.
\end{equation}
Then $\theta \in C_h$ if and only if
\begin{equation}
\label{cos_th}
\cos \theta = \frac{1}{2} \times \frac{|n^2 + 2mn - 2m^2|}{n^2-mn+m^2}
\end{equation}
for some $m,n \in \zed$ such that $\gcd(m,n)=1$, $1 \leq \frac{m}{n} \leq 2$, and $3 \nmid (m+n)$. For each $\theta \in C_h$, let us write $\cht$ for the corresponding similarity class. Then $\cht=\scgt$, where
\begin{equation}
\label{min_lattice}
\Gamma_{\theta} = \frac{1}{2} \begin{bmatrix} m+n&m-2n \\ (m-n)\sqrt{3}&m\sqrt{3} \end{bmatrix} \zed^2 \subseteq \Lambda_h,
\end{equation}
for the integers $m,n$ corresponding to $\theta$ as above, and for each $\Gamma \in \cht$, 
\begin{equation}
\label{min_lattice_det}
|\Gamma| \geq |\Gamma_{\theta}| = n^2-mn+m^2,\ |\Lambda_h : \Gamma| = \frac{\det(\Gamma)}{\det(\Lambda_h)} \geq |\Lambda_h : \Gamma_{\theta}| =  (2m-n)n.
\end{equation}
In fact,
\begin{equation}
\label{min_lattice_det_1}
\cht = \left\{ \sqrt{k} A \Gamma_\theta \subseteq \Lambda_h : k \in \zed_{> 0},\ A \in O_2(\real) \right\},
\end{equation}
where all the possible values of $k$ and the corresponding matrices $A$ are explicitly described in Lemma~\ref{sim_class_desc} below, and all possible values of $|\Lambda_h: \Gamma|$ for $\Gamma \in \WR(\Lambda_h)$ are described in Corollary \ref{index_values} below. Due to the properties \eqref{min_lattice_det} and \eqref{min_lattice_det_1}, we call $\Gamma_{\theta}$ a \textup{minimal sublattice} in its similarity class.
\end{thm}

\begin{rem} \label{ideal} Notice in particular that $\pi/3 \in C_h$ with the corresponding pair $(m,n)=(1,1)$, and $\ch(\pi/3) = \left< \Lambda_h \right>$. As is indicated in \cite{sloane}, ideal sublattices of $\Lambda_h$ are precisely those that are similar to $\Lambda_h$, hence the ideal sublattices form only one similarity class $\ch(\pi/3)$ in the infinite set $\swrh$ of similarity classes of WR sublattices of $\Lambda_h$, parameterized by $C_h$. On the other hand, $\pi/2 \notin C_h$, since
$$n^2 + 2mn - 2m^2 = (m+n)^2-3m^2 \neq 0,$$
as $3 \nmid (m+n)$. This is the complete opposite of the situation for WR sublattices of $\zed^2$ (identified with the ring of Gaussian integers), as studied in \cite{me:wr1}: the similarity class of ideal sublattices of $\zed^2$ corresponds to the value of the angle $\pi/2$, while there is no similarity class of WR sublattices of $\zed^2$ corresponding to $\pi/3$.
\end{rem}
\smallskip

We prove Theorem \ref{sim_par} in section~\ref{norm_f}. Our proof is based on a parameterization of lattices in question in terms of integral solution to a certain Diophantine equation, given by a ternary quadratic form. We produce such a parameterization of solutions for a family of integral ternary quadratic form equations in section~\ref{dioph_e} with the use of a simple geometric argument. In section~\ref{angle_f} we show how a particular such quadratic form can be used to parameterize similarity classes of WR lattices by means of looking at the corresponding values of the angle. However the parameterization of section~\ref{angle_f} does not necessarily  produce a minimal lattice for each similarity class. The main goal of section~\ref{norm_f} then is to go one step further and produce a description of similarity classes in terms of minimal lattices.

In section~\ref{optimize} we discuss the three optimization questions for WR sublattices of $\Lambda_h$, that are analogous to the questions considered in \cite{sloane}. These questions are concerned with counting the number, as well as maximizing the minimal norm and signal-to-noise ratio of WR sublattices of $\Lambda_h$ of a fixed index. Given a sublattice $\Gamma \in \WR(\Lambda_h)$, we can regard its nonzero points as transmitters which interfere with the transmitter at the origin, and then a standard measure of the {\it total interference} of $\Gamma$ is given by $E_{\Gamma}(2)$, where
\begin{equation}
\label{epstein_z}
E_{\Gamma}(s) = \sum_{\bx \in \Gamma \setminus \{ \bo \}} \frac{1}{\|\bx\|^{2s}}
\end{equation}
is the Epstein zeta-function of $\Gamma$, and the signal-to-noise ratio of $\Gamma$ is defined by
\begin{equation}
\label{SNR_def}
\SNR(\Gamma) = 10 \log_{10} \frac{1}{9E_{\Gamma}(2)},
\end{equation}
as in \cite{sloane}. To maximize $\SNR(\Gamma)$ on the set of all WR sublattices of $\Lambda_h$ of a fixed index $J$ is the same as to minimize $E_{\Gamma}(2)$. In particular, we show (Lemma~\ref{SNR_min} below) that $\SNR(\Gamma)$ is maximized by the same sublattice of fixed index $J$ that maximizes $|\Gamma|$, and vice versa. This is not always so for non-WR sublattices of $\Lambda_h$, as demonstrated in \cite{sloane}.

Finally, in section~\ref{combin} we discuss a combinatorial structure on the set of all similarity classes of WR sublattices of $\Lambda_h$, induced by the action of a certain submonoid of $\GL_3(\zed)$.
\bigskip

\section{A Diophantine equation}
\label{dioph_e}

In this section we use a simple geometric idea to construct an explicit parameterization of integral zeros of a certain integral ternary quadratic form. We later use this parameterization to prove Theorem \ref{sim_par}.
    
\begin{lem} \label{dioph} Consider the Diophantine equation
\begin{equation}
\label{d1}
\alpha x^2 + \beta xy + \gamma y^2 = \delta z^2,
\end{equation}
where $\alpha,\beta,\gamma,\delta \in \ZZ$ with $\beta^2 \neq 4\alpha\gamma$ and $\delta \neq 0$. Then either this equation has no integral solutions with $z \neq 0$, or all such solutions $(x,y,z)$ of \eqref{d1} are rational multiples of
\begin{equation}
\label{sol}
\begin{split}
x & = \gamma n (an - 2bm) - (\alpha a + \beta b) m^2, \\
y & = \alpha m (bm - 2an) - (\gamma b + \beta a) n^2, \\ 
z & = \pm c (\alpha m^2 + \beta mn + \gamma n^2),
\end{split}
\end{equation}
where $m,n \in \zed$ with $\gcd(m,n)=1$ and $m \geq 0$; here $(a,b,c)$ is any integral solution to \eqref{d1} with $c \neq 0$. In this later case, every multiple of \eqref{sol} is a solution to \eqref{d1} by homogeneity of the equation \eqref{d1}.
\end{lem}
    
\proof
Suppose there exists an integer solution $(a,b,c)$ to \eqref{d1} with $c \neq 0$, and consider the rational curve
\begin{equation}
\label{e1}
\frac{\alpha}{\delta}\, u^2 + \frac{\beta}{\delta}\, uv + \frac{\gamma}{\delta}\, v^2 = 1,
\end{equation}
where $u=x/z$, $v=y/z$. The point $(a/c, b/c)$ lies on this curve, and let $(u,v)$ be any other rational point on the curve. If $v = b/c$, then we must have $\alpha \neq 0$ and $u = -a/c - b\beta/c\alpha$, and so the point $(u,v)$ corresponds to the solution $(-a\alpha-b\beta,b\alpha,c\alpha)$ of \eqref{d1}, obtained from \eqref{sol} when $m=1,\ n=0$. Otherwise, there exists a unique line with rational slope through the points $(a/c, b/c)$ and $(u,v)$. Therefore there exists $m/n \in \que$ such that
\begin{equation}
\label{u_line}
u = \frac{m}{n} \left( v - \frac{b}{c} \right) + \frac{a}{c}.
\end{equation}
Substituting \eqref{u_line} into \eqref{e1} and solving for $v$, we obtain:
\begin{equation}
\label{v}
v = \frac{ \alpha m (bm - 2an) - (\gamma b + \beta a) n^2 }{ c (\alpha m^2 + \beta mn + \gamma n^2) }.
\end{equation}
Now substituting \eqref{v} into \eqref{u_line}, we have:
\begin{equation}
\label{u}
u = \frac{ \gamma n (an - 2bm) - (\alpha a + \beta b) m^2 }{ c (\alpha m^2 + \beta mn + \gamma n^2) }.
\end{equation}

        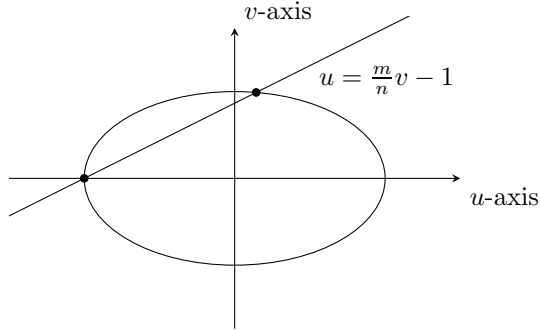
\begin{figure}[h]
            \centering
            \begin{tikzpicture}[>=stealth]
                % U & V Axes %
                \draw[->] (-3, 0) -- (3, 0) node[below right] {$u$-axis};
                \draw[->] (0, -2) -- (0, 2) node[above right] {$v$-axis};
                
                % Ellipse %
                \draw (0, 0) ellipse (2 and 1.1547005381);
                
                % Line %
                \begin{scope}
                    \clip (-3, -2) rectangle (3, 2.15);
                    \draw (-4.2857142857, -1.1428571429) -- (2.5714285714, 2.2857142857);
                \end{scope}
                
                % Intersections %
                \filldraw (-2, 0) circle (.05);
                \filldraw (0.2857142857, 1.1428571429) circle (.05);
                
                % Line Label %
                \draw (1, 1) node[above right] {$u =\frac{m}{n}v - 1$};
            \end{tikzpicture}
            \caption{Algorithm for $x^2 + 3y^2 = z^2$ with the solution $(-1,0,1)$.} \label{fig_a}
        \end{figure}

\noindent
As $m,n$ range over all coprime pairs of integers, \eqref{u} and \eqref{v} give coordinates of all rational points $(u,v)$ on the curve \eqref{e1}, and therefore every integral solution $(x,y,z)$ of \eqref{d1} is a rational multiple of \eqref{sol} for some $m,n$.
\endproof
    
\begin{rem} \label{proj} Notice that for every two different relatively prime pairs of integers $m_1,n_1$ and $m_2,n_2$, the corresponding solutions $(x_1,y_1,z_1)$ and $(x_2,y_2,z_2)$ of \eqref{d1} as in \eqref{sol} represent different projective points, i.e., there does not exist a real number $t$ such that $(x_1,y_1,z_1) = t (x_2,y_2,z_2)$. Hence \eqref{sol} describes all the distinct rational projective points on the hypersurface given by the equation \eqref{d1}. We will be especially interested in {\it primitive} integral solutions to equations of the form  \eqref{d1}, i.e., representatives $(x,y,z) \in \zed^3$ of their projective points described by \eqref{sol} with $\gcd(x,y,z) = 1$. 
\end{rem}
\bigskip

\section{The angle form}
\label{angle_f}

In this section we construct our first parameterization of WR sublattices of $\Lambda_h$.

\begin{lem} \label{angle_ch} For every $\Gamma \in \WR(\Lambda_h)$,
\begin{equation}
\label{cos_sin}
\cos \theta(\Gamma) = \frac{p}{q} \leq \frac{1}{2}, \quad \sin \theta(\Gamma) = \frac{r}{q}\sqrt{3} \geq \frac{\sqrt{3}}{2},
\end{equation}
where the triple $(p,r,q) \in \zed^3_{\geq 0}$ is relatively prime.
\end{lem}

\proof
Let $\Gamma$ be a WR sublattice of $\Lambda_h$. Then there exist $a,b,c,d \in \zed$ such that
 \[
        \Gamma = 
            \begin{bmatrix}
                1 & -\frac{1}{2} \\
                0 & \frac{\sqrt{3}}{2}
            \end{bmatrix}
            \begin{bmatrix}
                a & c \\
                b & d
            \end{bmatrix} \zed^2 = \begin{bmatrix} a-b/2 & c-d/2 \\ b\sqrt{3}/2 & d\sqrt{3}/2 \end{bmatrix} \zed^2,
    \]
where
$$\bx := \begin{bmatrix} a-b/2 \\ b\sqrt{3}/2 \end{bmatrix}, \quad \bwy := \begin{bmatrix} c-d/2 \\   d\sqrt{3}/2 \end{bmatrix}$$
is a minimal basis for $\Gamma$. Hence 
$$\|\bx\|^2 = a^2-ab+b^2 = c^2-cd+d^2 = \|\bwy\|^2,$$
and
$$\frac{1}{2} \geq \cos \theta(\Gamma) = \frac{\bx^t \bwy}{\|\bx\| \|\bwy\|} = \frac{a (2c-d) + b (2d-c)}{2(a^2 - ab + b^2)} \in \que.$$
Then
\begin{eqnarray*}
\sin \theta(\Gamma) & = & \sqrt{1 - \cos^2 \theta(\Gamma)} =  \sqrt{ 1 -  \frac{\left( 2ac + 2bd - ad - bc \right)^2}{4 (a^2 - ab + b^2)^2} } \\
& = & \frac{ad - bc}{2(a^2 - ab + b^2)} \sqrt{3} \geq \frac{\sqrt{3}}{2}.
\end{eqnarray*}
Therefore \eqref{cos_sin} is satisfied with
\begin{equation}
\label{pqr}
p = a (2c-d) + b (2d-c), \quad r = ad - bc, \quad q = 2 (a^2 - ab + b^2),
\end{equation}
which completes the proof of the lemma.
\endproof
\smallskip

Notice, on the other hand, that if the triple $(p,r,q)$ is as in \eqref{cos_sin} for some $\Gamma \in \WR(\Lambda_h)$, then
\begin{equation}
\label{angle_form}
p^2+3r^2=q^2,
\end{equation}
which is a particular instance of the equation \eqref{d1} with $\alpha = \delta = 1$, $\beta=0$, and $\gamma=3$. We will refer to the binary integral quadratic form $p^2+3r^2$ on the left hand side of \eqref{angle_form} as the {\it angle form} because of the above connection of the integral solutions of \eqref{angle_form} to values of trigonometric functions of the angles of lattices from $\WR(\Lambda_h)$. In fact, Lemma \ref{dioph} (see also Figure \ref{fig_a} above) guarantees that all projectively distinct integer solutions $(p,r,q) \in \ZZ^3_{>0}$ of \eqref{angle_form} are given by
\begin{equation}
\label{mn}
p = m^2 - 3n^2, \quad r = 2mn, \quad q = m^2 + 3n^2,
\end{equation}
where $m$ and $n$ are coprime integers. If we only consider triples $(p,r,q)$ as in \eqref{mn} with $\sqrt{3} < m/n \leq 3$, then there exists an angle $\theta \in [\pi/3, \pi/2)$ such that
\begin{equation}
\label{cos_sin_1}
\cos \theta = \frac{p}{q}, \quad \sin \theta = \frac{r}{q}\sqrt{3},
\end{equation}
and vice versa. Then $\theta \in C_h$ defines the similarity class $\cht$ of WR sublattices of $\Lambda_h$, as described in the statement of Theorem \ref{sim_par}. We will now construct a sublattice $\gat$ of $\Lambda_h$ in $\cht$ using our parameterization.
    
\begin{lem} \label{Gamma} Let $m$ and $n$ be coprime positive integers with $\sqrt{3} < \frac{m}{n} \leq 3$, and let $(p,r,q)$ be as in \eqref{mn}. Let $\theta$ be as in \eqref{cos_sin_1}. Then the lattice
\begin{equation}
\label{angle_lattice}
\gat := \begin{bmatrix}
                    m & m \\
                    n\sqrt{3} & -n\sqrt{3}
                \end{bmatrix} \zed^2
        \end{equation}
is in $WR(\Lambda_h)$, and has the following properties:
        \[
            \abs{\gat} = q, \quad
            \det (\gat) = r\sqrt{3}, \quad
            |\Lambda_h : \gat| = 2r, \quad \text{and} \quad
            \gat \in \cht.
        \]
Hence $\cht = \left< \gat \right>$.
\end{lem}
    
\proof 
First notice that $\gat$, as defined in \eqref{angle_lattice}, is given by
        \[
            \gat =
                \begin{bmatrix}
                    1 & -\frac{1}{2} \\
                    0 & \frac{\sqrt{3}}{2}
                \end{bmatrix}
                \begin{bmatrix}
                    m+n & m-n \\
                    2n & -2n
                \end{bmatrix} \zed^2,
        \]
        hence it is a sublattice of $\Lambda_h$. Also
        \[
            \det (\gat) = 2mn\sqrt{3} = r\sqrt{3},
        \]
        and
        \[
            |\Lambda_h : \gat| = \frac{\det\Gamma}{\det\Lambda_h} = 4mn = 2r.
        \]
        Now let
        \[
            \bx =
                \begin{bmatrix}
                    m \\
                    n\sqrt{3}
                \end{bmatrix},
            \quad
            \bwy =
                \begin{bmatrix}
                    m \\
                    -n\sqrt{3}
                \end{bmatrix},
        \]
and notice that $\|\bx\| = \|\bwy\| = \sqrt{m^2+3n^2} = \sqrt{q}$. Moreover, if $\nu$ is the angle between these two vectors, then
        \[
            \cos \nu = \frac{ \bx^t \bwy }{ \|\bx\| \|\bwy\| } = 
                \frac{ m^2 - 3n^2}{q} =
                \frac{p}{q},
        \]
hence this angle is precisely $\theta$. The condition $\sqrt{3} < \frac{m}{n} \leq 3$ implies that $\pi/3 \leq \theta < \pi/2$, meaning that for any $s,t \in \zed$
$$\|s \bx + t \bwy \|^2 = s^2 \|\bx\|^2 + 2st \|\bx\| \|\bwy\| \cos \theta + t^2 \|\bwy\|^2 \geq \|\bx\|^2 = \|\bwy\|^2,$$
hence $\bx,\bwy$ form a minimal basis for $\gat$. Therefore $\abs{\gat} = q$, and $\gat \in \cht$.
\endproof
\smallskip
    
Thus similarity classes of WR sublattices of $\Lambda_h$ are in bijective correspondence with triples $(p,r,q)$ as defined in \eqref{mn} satisfying $\sqrt{3} < m/n \leq 3$ and $\gcd(m,n)=1$. In other words, we obtained a parameterization of the set $C_h$ by a subset of projectively distinct zeros of the angle form. The construction of Lemma \ref{Gamma} however does not always produce {\it minimal} WR sublattices of $\Lambda_h$ in the sense of \eqref{min_lattice_det} of Theorem \ref{sim_par}. For instance,  $\Lambda_h \in \ch(\pi/3)$, while 
\[
\Omega_{\pi/3} = \begin{bmatrix}
                    3 & 3 \\
                    \sqrt{3} & -\sqrt{3}
                \end{bmatrix} \zed^2 
\]
with $|\Omega_{\pi/3}| = 12$ (recall that $|\Lambda_h| = 1$). Our next goal is to provide a parameterization of $\swrh$ by minimal sublattices, which we do in the next section.
\bigskip

\section{The norm form}
\label{norm_f}

The main goal of this section is to prove Theorem \ref{sim_par}. We start by considering another particular instance of the equation \eqref{d1} with $\alpha = \gamma = \delta = 1$ and $\beta = -1$, which we write as
\begin{equation}
\label{nf}
a^2 - ab + b^2 = c^2,
\end{equation}
where the left hand side of \eqref{nf} is the {\it norm form} of the hexagonal lattice with respect to the basis matrix $\begin{bmatrix} 1&-\frac{1}{2} \\ 0&\frac{\sqrt{3}}{2} \end{bmatrix}$, which is precisely the norm in the ring of Eisenstein integers $\zed[\omega]$. We will call every solution $(a,b,c) \in \zed^3_{\geq 0} \setminus \{(0,0,0)\}$ to \eqref{nf} an {\it Eisenstein triple}, and we call an Eisenstein triple $(a,b,c)$ {\it primitive} if $a \leq b$ and $\gcd(a,b,c) = 1$. Then we have the following simple corollary of Lemma \ref{dioph}.

\begin{cor} \label{spec_param} All Eisenstein triples are positive rational multiples of
\begin{equation}
\label{sol2}
a = m (2n - m), \quad b = n (2m - n), \quad c = m^2 - mn + n^2,
\end{equation}
where $m,n \in \zed_{> 0}$ with $\gcd(m,n) = 1$ and $\frac{1}{2} \leq \frac{m}{n} \leq 2$. Moreover, every integer triple $(a,b,c)$ that is a positive rational multiple of \eqref{sol2} for some $m,n$ satisfying the above conditions is an Eisenstein triple by homogeneity of the equation \eqref{nf}.
\end{cor}

\proof
We apply Lemma \ref{dioph} to the equation \eqref{nf} with the particular solution $(-1,-1,1)$ as in Figure \ref{fig_b} below.
        \begin{figure}[h]
            \centering
            \begin{tikzpicture}[>=stealth]
                % U & V Axes %
                \draw[->] (-3, 0) -- (3, 0) node[below right] {$u$-axis};
                \draw[->] (0, -3) -- (0, 3) node[above right] {$v$-axis};

                % Ellipse %
                \draw[rotate=45] (0, 0) ellipse (2.82842712474619 and 1.632993161855452);

                % Line %
                \begin{scope}
                    \clip (-3, -3) rectangle (3, 3.15);
                    \draw (-4.8571428571, -6.2857142857) -- (3.7142857143, 6.5714285714);
                \end{scope}

                % Intersections %
                \filldraw (-2, -2) circle (.05);
                \filldraw (0.8571428571, 2.2857142857) circle (.05);

                % Line Label %
                \draw (1.25, 2.45) node[above right] {$u =\frac{m}{n}(v+1) - 1$};
            \end{tikzpicture}
            \caption{Algorithm for $x^2 - xy + y^2 = z^2$ with the solution $(-1,1,1)$.} \label{fig_b}
        \end{figure}
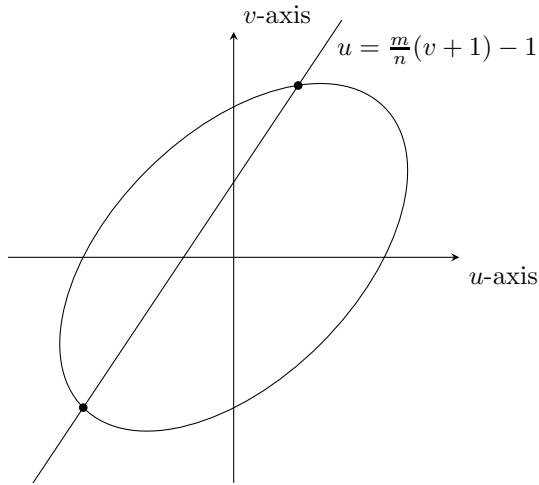
To obtain the restrictions on $m$ and $n$, first notice that since $(a,b,c) \in \zed^3_{\geq 0}$ and $m \geq 0$ we have $0 \leq 2n - m$, hence $n \geq \frac{m}{2} \geq 0$ and $\frac{m}{n} \leq 2$ ($n > 0$ since $n=0$ would imply $m=n=0\ \Rightarrow\ a=b=c=0$). Since $b \geq 0$, we must have $0 \leq 2m - n$. Since $n>0$, we have $m>0$, and hence $\frac{1}{2} \leq \frac{m}{n}$.
\endproof

Here is a brief outline of our strategy for the remainder of this section. We will first use Corollary \ref{spec_param} to give a parameterization of the primitive Eisenstein triples. We will then relate primitive Eisenstein triples to solutions of the angle form equation, thus obtaining a new parameterization for the similarity classes of WR sublattices of $\Lambda_h$. This new parameterization will then be used to produce minimal lattices for the similarity classes.
\smallskip

Since we are trying to produce primitive Eisenstein triples, we may assume that integers $m,n$ in \eqref{sol2} are relatively prime. We start by proving that in this case $\gcd(a,b,c) = 1$ or $3$. 

\begin{lem} \label{coprimality} Let $(a,b,c)$ be a triple of the form \eqref{sol2} with relatively prime $m,n \in \zed$, then $\gcd(a,b,c) = 1$ or $3$. More precisely:
\begin{enumerate}
\item $\gcd(a, b, c) = 3$ if and only if $3 \mid (m + n)$.
\item $\gcd(a, b, c) = 1$ if and only if $3 \nmid (m + n)$.
\end{enumerate}
\end{lem}

\proof
Suppose first that $m + n = 3k$ for some integer $k$, then
$$a = 3m (2k - m), \quad b = 3n (2k - n), \quad c = 3 (km - mn + kn),$$
so $3 \mid \gcd(a,b,c)$. Now suppose that there is some prime $p$ dividing $\frac{\gcd(a,b,c)}{3}$, then $p \mid m(2k-m), n(2k-n)$. Since $\gcd(m,n) = 1$, it must be that either $p \mid m, 2k-n$ or $p \mid 2k-m, n$. Notice that $n = 3k-m$ and $2k-n = m-k$. First suppose that $p \mid m, 2k-n$, then
$$p \mid m-(m-k) = k\ \Longrightarrow\ p \mid 2k-(2k-n) = n,$$
which is a contradiction. Next suppose that $p \mid 2k-m, n$, then
$$p \mid (3k-m)-(2k-m) = k\ \Longrightarrow\ p \mid 2k-(2k-m)=m,$$
which is a contradiction. Therefore $\gcd(a,b,c) = 3$.

Conversely, if $\gcd(a, b, c) = 3$, then
$$3 \mid a+b = -((m + n)^2 - 6mn),$$
hence $3 \mid m+n$.

Next assume that $3 \nmid (m + n)$. Suppose that $\gcd(a, b, c) > 1$, then there exists a prime $p$ such that $p \mid a, b$. Suppose that $p \mid (m-n)$, then
$$p \mid (b + (m-n)^2) = m^2,\ \ p \mid (a + (m-n)^2) = n^2,$$
which contradicts co-primality of $m,n$. Therefore $p \nmid (m-n)$. Moreover, $p \nmid m,n$: if, for instance, $p \mid m$, then $p \mid (2nm-b) = n^2$, which is a contradiction; similarly, if $p \mid n$, then $p$ must divide $m$. Now notice that 
$$p \mid (b-a) = (m-n) (m+n),$$ 
and since $p \nmid (m-n)$, $p$ must divide $m+n$. Therefore 
$$p \mid (b+(m+n)^2)=m(m+4n),\ \ p \mid (a + (m+n)^2) = n (n + 4m).$$
Since $p \nmid m,n$, we find that $p \mid (m+4n), (n+4m)$, and hence 
$$p \mid ((n+4m) - (m+4n)) = 3(m-n).$$
Since $p \nmid (m-n)$, it must be the case that $p = 3$, contradicting our assumption that $3 \nmid (m+n)$. Thus $\gcd(a,b,c) = 1$.

Conversely, suppose $\gcd(a,b,c) = 1$. Suppose $3 \mid (m+n)$. Then
$$3 \mid (m-n) (m+n) = b-a,\ \ 3 \mid -((m + n)^2 - 6mn) = a+b,$$
and so
$$3 \mid (b-a)+(b+a) = 2b,\ \ 3 \mid (b-a)-(b+a) = -2a.$$
Hence $3 \mid a,b$, which is a contradiction, so $3 \nmid (m+n)$. This completes the proof of the lemma.
\endproof

\begin{rem} \label{rem_barnes} A parameterization of primitive Eisenstein triples partially similar to our Corollary \ref{spec_param} and Lemma \ref{coprimality} (but using different arguments) has been obtained in \cite{barnes}. We include our results here since the details of this parameterization are important to our description of similarity classes of WR sublattices of $\Lambda_h$.
\end{rem}
\smallskip

We can now give a parameterization of primitive Eisenstein triples. Notice that for each primitive Eisenstein triple $(a,b,c)$, $(b-a,b,c)$ is also a primitive Eisenstein triple, and $(b-(b-a),b,c) = (a,b,c)$, i.e., the map taking $(a,b,c)$ to $(b-a,b,c)$ is an involution on the set of primitive Eisenstein triples. We will call two triples like this {\it associated}, and write $\left< a,b,c \right>$ for the associated pair. Then the set of all primitive Eisenstein triples can be split into a collection of associated pairs.

\begin{lem} \label{norm_param} A pair of vectors $(a,b,c), (b-a,b,c) \in \zed^3_{\geq 0}$ is an associated pair of primitive Eisenstein triples if and only if precisely one of these triples is given by \eqref{sol2} for some integers $m,n$ satisfying
\begin{enumerate}
\item $m,n > 0$ and $\gcd(m,n) = 1$
\item $1 \leq \frac{m}{n} \leq 2$
\item $3 \nmid (m + n)$.
\end{enumerate}
\end{lem}

\proof
Suppose first that $(a,b,c)$ is as in \eqref{sol2} with $m,n$ satisfying the above-stated conditions (1), (2), (3). Then $(a,b,c)$ is an Eisenstein triple by Corollary \ref{spec_param}, and $\gcd(a,b,c) = 1$ by Lemma \ref{coprimality}. Moreover,
$$a = 2mn - m^2 \leq 2mn - n^2 = b,$$
and hence $(a,b,c)$ is a primitive Eisenstein triple.

Now assume that $(a,b,c)$ is a primitive Eisenstein triple, then so is $(b-a,b,c)$. Corollary \ref{spec_param} implies that
$$(a,b,c) = \frac{1}{g} (a',b',c'),$$
where $(a',b',c')$ is an Eisenstein triple  given by \eqref{sol2}  with $m,n$ satisfying condition (1) and $\frac{1}{2} \leq \frac{m}{n} \leq 2$, and $g=\gcd(a',b',c')$. First suppose that $\frac{m}{n} < 1$. This means that
\begin{equation}
\label{ab_ineq}
a = 2mn - m^2 > 2mn - n^2 = b,
\end{equation}
which contradicts $(a,b,c)$ being primitive. Therefore the integers $m,n$ must also satisfy condition (2). Now Lemma \ref{coprimality} implies that $m,n$ satisfy condition (3) if and only if $g=1$, and $g=3$ otherwise. Suppose $g=3$. Then it can be easily verified that $(a_1,b_1,c_1):=(b-a,b,c)$ is given by \eqref{sol2} with positive integers
$$m_1 = \frac{m+n}{3},\ n_1 = \frac{2m-n}{3}.$$
Suppose that an integer $t \mid m_1,n_1$. Then
$$t \mid m_1+n_1 = m, 2m_1-2n_1=n,$$
but $\gcd(m,n)=1 \Rightarrow \gcd(m_1,n_1)=1$, hence $m_1,n_1$ satisfy condition (1). By the same argument as in \eqref{ab_ineq} with $a_1,b_1,m_1,n_1$ instead of $a,b,m,n$ we conclude that $m_1,n_1$ satisfy condition (2). Finally, suppose that $3 \mid (m_1+n_1) = m$. Since $3 \mid (m+n)$, we conclude that $3 \mid n$, meaning that $m,n$ are not relatively prime, which is a contradiction. Therefore $m_1,n_1$ satisfy condition (3) as well.

Now suppose that $g=1$, then it is easy to verify that
$$(b-a,b,c) = \frac{1}{3} (a_2,b_2,c_2),$$
where $(a_2,b_2,c_2)$ is given by \eqref{sol2} with positive integers
$$m_2=m+n,\ n_2 = 2m-n.$$
By the argument in the proof of Lemma \ref{dioph} (see also Remark \ref{proj}), there cannot exist another pair $m_2',n_2'$ parameterizing the triple $(b-a,b,c)$. Hence we showed that {\it precisely} one of the triples in an associated pair is given by \eqref{sol2} with integers $m,n$ satisfying (1)-(3). This completes the proof of the lemma.
\endproof

Our next step is to relate primitive Eisenstein triples to solutions of the angle form equation \eqref{angle_form}. For each vector $(x,y,z) \in \real^3$, we will write $[x,y,z]$ to denote the corresponding projective point. First define
\begin{eqnarray}
\label{angle_set}
\AA & = & \{[p,r,q] :  (p,r,q) \text{ satisfying \eqref{mn}, } \sqrt{3} < m/n \leq 3 \} \nonumber \\
& = & \{ [p,r,q] : (p,r,q) \text{ satisfying \eqref{angle_form}, } 0 < p/q \leq 1/2 \}.
\end{eqnarray}
Next, for each associated pair $\left< a,b,c \right>$ of primitive Eisenstein triples, let us write $\{a,b,c\}$ for the associated pair of the corresponding projective points, $[a,b,c]$ and $[b-a,b,c]$, and define $\EE$ to be the set of all such associated pairs. Let $T : \real^3 \to \real^3$ be a bijective linear map, given by the matrix
$$\begin{bmatrix} -2&1&0 \\ 0&1&0 \\ 0&0&2 \end{bmatrix}.$$
Notice that for each $\{a,b,c\} \in \EE$ with $\gcd(a,b,c)=1$
$$T(a,b,c) = (b-2a,b,2c), \quad T(b-a,b,c) = (-(b-2a),b,2c),$$
and
$$(b-2a)^2+3b^2=(2c)^2,$$
meaning that precisely one triple from the associated pair $\{a,b,c\}$ maps to the triple $(|b-2a|,b,2c)$ under $T$. Since $a \leq b$ and $\gcd(a,b)=1$ (meaning in particular that $2a \neq b$),
$$0 < |b-2a|^2 = 4a (a - b) + b^2 \leq a (a - b) + b^2 = c^2,$$
therefore $0 < |b-2a|/2c \leq 1/2$, and so $[|b-2a|,b,2c] \in \AA$. Now consider the map $T^{-1} : \real^3 \to \real^3$, given by the inverse matrix of $T$
$$\begin{bmatrix} -1/2&1/2&0 \\ 0&1&0 \\ 0&0&1/2 \end{bmatrix},$$
and notice that for every $[p,r,q] \in \AA$,
$$T^{-1}(p,r,q) = \left( \frac{r-p}{2}, r, \frac{q}{2} \right)$$
is a solution to \eqref{nf}. Moreover, if
$$p = m^2 - 3n^2, \quad r = 2mn, \quad q = m^2 + 3n^2,$$
where $m$ and $n$ are coprime integers with $m = \alpha n$ for some real $\sqrt{3} < \alpha \leq 3$, then
$$r-p = (2\alpha - \alpha^2 + 3)n^2 \geq 0,$$
and $(r-p)/2 \leq r$. Hence $T^{-1}(p,r,q)$ is a representative of a projective point defined by some primitive Eisenstein triple. We have proved the following lemma.

\begin{lem} \label{sets_biject} There is a bijective correspondence between the sets $\AA$ and $\EE$ as described above.
\end{lem}

Now this bijection can be easily used to relate primitive Eisenstein triples to well-rounded sublattices of the hexagonal lattice.

\begin{cor} \label{norm_angle_cor} Let $C_h$ be as in \eqref{set_ch}, and for every $\theta \in C_h$ let $\cht$ be the corresponding similarity class of WR sublattices of $\Lambda_h$, as in the statement of Theorem \ref{sim_par}. Then for each $\cht$ there exist exactly two primitive Eisenstein triples $(a,b,c)$ and $(b-a,b,c)$ such that
\begin{equation}
\label{norm_angle}
\cos\theta = \frac{\abs{b - 2a}}{2c} = \frac{\abs{b - 2 (b-a)}}{2c}.
\end{equation}
Conversely, for each primitive Eisenstein triple $(a,b,c)$ there exists a unique similarity class $\cht$ with $\theta \in C_h$ satisfying \eqref{norm_angle}.
\end{cor}

\proof
Notice that the quotients in \eqref{norm_angle} only depend on the associated pair of projective points $\{a,b,c\}$, not on actual choice of representatives of these points. The statement of the corollary now follows by combining Lemma \ref{sets_biject} with the parameterization of section~\ref{angle_f}.
\endproof

Next we use the parameterization given in Lemma \ref{norm_param} to construct a WR sublattice in each similarity class $\cht$ with $\theta$ as in \eqref{norm_angle}.

\begin{lem} \label{minimal_gamma} Let $(a,b,c)$ be a primitive Eisenstein triple, and let $m,n$ be the integers parameterizing either $(a,b,c)$ or $(b-a,b,c)$, as defined in Lemma \ref{norm_param}. Let $\theta$ be as given by \eqref{norm_angle}. Define
\begin{equation}
\label{minimal_sub}
\Gamma_\theta = \begin{bmatrix} 1 & -\frac{1}{2} \\ 0 & \frac{\sqrt{3}}{2} \end{bmatrix}
\begin{bmatrix} m & m-n \\ m-n & m \end{bmatrix} \zed^2 = 
\frac{1}{2} \begin{bmatrix} m+n&m-2n \\ (m-n)\sqrt{3}&m\sqrt{3} \end{bmatrix} \zed^2.
\end{equation}
Then $\Gamma_{\theta} \in \WR(\Lambda_h)$ is such that $\cht = \left< \Gamma_{\theta} \right>$. Moreover,
\[
\abs{\Gamma_\theta} = c, \quad
\det \Gamma_\theta = b\frac{\sqrt{3}}{2}, \quad
|\Lambda_h : \Gamma_{\theta}| = b.
\]
\end{lem}

\proof
By definition \eqref{minimal_sub}, $\Gamma_{\theta}$ is a sublattice of $\Lambda_h$. First we see that
$$\det \Gamma_\theta = \det \begin{bmatrix} m & m-n \\ m-n & m \end{bmatrix} \det \Lambda_h = n (2n - m) \frac{\sqrt{3}}{2} = b \frac{\sqrt{3}}{2},$$
and so
$$|\Lambda_h : \Gamma_\theta| = \frac{\det \Gamma_{\theta}}{\det \Lambda_h} = b.$$
Now let
$$\bx = \frac{1}{2} \begin{bmatrix} m+n \\ (m-n) \sqrt{3} \end{bmatrix}, \quad \bwy = \pm \frac{1}{2} \begin{bmatrix} m - 2n \\ m \sqrt{3} \end{bmatrix},$$
where the $\pm$ choice in the definition of $\bwy$ is such that $\bx^t \bwy>0$. Notice that $\|\bx\|^2 = \|\bwy\|^2 = m^2 - mn + n^2 = c$. Moreover, the cosine of the angle between these two vectors is equal to
$$\frac{\bx^t \bwy}{\|\bx\| \|\bwy\|} = \pm \frac{2m^2 - 2mn - n^2}{2c} = \frac{|b - 2a|}{2c},$$
hence this angle is precisely $\theta$. By construction, $\pi/3 \leq \theta < \pi/2$, which means that any integral linear combination of the vectors $\bx$ and $\bwy$ has norm at least as large as the norm of these vectors, hence they form a minimal basis for $\Gamma_{\theta}$. This implies that $|\Gamma_{\theta}| = c$ and $\Gamma_{\theta} \in \cht$, which means that $\cht = \left< \Gamma_{\theta} \right>$.
\endproof

Now we show that the sublattice constructed in Lemma \ref{minimal_gamma} is minimal in its similarity class, as defined in the statement Theorem \ref{sim_par}.

\begin{lem} \label{gstm} Let $\Gamma_\theta$ be defined as in \eqref{minimal_sub} and let $\Omega \in \cht$, then $|\Gamma_{\theta}| \leq |\Omega|$. Hence $\Gamma_{\theta}$ is a minimal sublattice in the similarity class $\cht$.
\end{lem}

\proof
Since $\Omega \sim \Gamma_{\theta}$, there exists $0 \neq \alpha \in \real$ and a $2 \times 2$ real orthogonal matrix $A$ such that
$$\Omega = \alpha A \Gamma_\theta.$$
Notice that $\alpha^2 c = |\Omega| \in \zed$ and $\alpha^2 b = |\Lambda_h : \Omega| \in \zed$, hence $\alpha^2 \in \que$; let us write $\alpha^2 = p/q$ with $\gcd(p,q) = 1$. Then $q^2 \mid b$ and $q^2 \mid c$, but $\gcd(b,c) = 1$ and so $\alpha^2 \in \zed_{\geq 0}$. Therefore $|\Gamma_{\theta}| \leq \alpha^2 |\Gamma_{\theta}| = |\Omega|$.
\endproof

Therefore for each $\theta \in C_h$,
\begin{equation}
\label{cht_0}
\cht = \left\{ \sqrt{k} A \Gamma_\theta \subseteq \Lambda_h : k \in \zed_{> 0},\ A \in O_2(\real) \right\}.
\end{equation}
What can be said about $k$ and $A$? We have the following lemma.

\begin{lem} \label{sim_class_desc} Let
$$D =  \{1\} \cup \{ d = p_1 \dots p_s : p_1,\dots,p_s \text{ distinct primes}  \equiv 1 (\md 3) \},$$
and for each $d \in D$ let 
$$S(d) = \{ (p,r,q) : p^2+3r^2 = dq^2 \}$$
as given by the parameterization of Lemma \ref{dioph} with $\alpha=1,\beta=0,\gamma=3,\delta=d$. A lattice $\Gamma \in \cht$ if and only if $\Gamma = \sqrt{k} A \Gamma_\theta$ where one of the following two conditions hold:
\begin{enumerate}
\item $k = j^2 d$ for some $j \in \zed_{>0}$, $d \in D$, and
$$A = \begin{bmatrix} \frac{p}{q \sqrt{d}} & -\frac{r  \sqrt{3}}{q \sqrt{d}} \\ \frac{r  \sqrt{3}}{q \sqrt{d}} & \frac{p}{q \sqrt{d}} \end{bmatrix} \text{ or } \begin{bmatrix} \frac{p}{q \sqrt{d}} & \frac{r  \sqrt{3}}{q \sqrt{d}} \\ \frac{r  \sqrt{3}}{q \sqrt{d}} & -\frac{p}{q \sqrt{d}} \end{bmatrix}$$
for some $(p,r,q) \in S(d)$ with $q \mid j$.
\item $k = 3 j^2 d$ for some $j \in \zed_{>0}$, $d \in D$, and
$$A = \begin{bmatrix} \frac{r  \sqrt{3}}{q \sqrt{d}} & -\frac{p}{q \sqrt{d}} \\ \frac{p}{q \sqrt{d}} & \frac{r  \sqrt{3}}{q \sqrt{d}} \end{bmatrix} \text{ or } \begin{bmatrix} \frac{r  \sqrt{3}}{q \sqrt{d}} & \frac{p}{q \sqrt{d}} \\ \frac{p}{q \sqrt{d}} & -\frac{r  \sqrt{3}}{q \sqrt{d}} \end{bmatrix}$$
for some $(p,r,q) \in S(d)$ with $q \mid j$.
\end{enumerate}
\end{lem}

\proof
Suppose $\Gamma \in \cht$, so $\Gamma = \sqrt{k} A \Gamma_\theta$ for some $k \in \zed_{> 0}$ and  $A \in O_2(\real)$.  First notice that
$$A = \begin{bmatrix} \cos t & -\sin t \\ \sin t & \cos t \end{bmatrix} \text{ or } \begin{bmatrix} \cos t & \sin t \\ \sin t & -\cos t \end{bmatrix} \text{ for some } 0 \leq t < 2\pi.$$
Therefore, either
$$\Gamma = \frac{\sqrt{k}}{2} \begin{bmatrix}  (m+n) \cos t - (m-n) \sqrt{3} \sin t & (m-2n) \cos t - m \sqrt{3} \sin t \\ (m+n) \sin t + (m-n) \sqrt{3} \cos t & (m-2n) \sin t + m \sqrt{3} \cos t \end{bmatrix} \zed^2,$$
or
$$\Gamma = \frac{\sqrt{k}}{2} \begin{bmatrix}  (m+n) \cos t + (m-n) \sqrt{3} \sin t & (m-2n) \cos t + m \sqrt{3} \sin t \\ (m+n) \sin t - (m-n) \sqrt{3} \cos t & (m-2n) \sin t - m \sqrt{3} \cos t \end{bmatrix} \zed^2,$$
where $m,n$ are as in Lemma \ref{minimal_gamma}; in any case, $\Gamma$ is a sublattice of $\Lambda_h = \begin{bmatrix} 1&-\frac{1}{2} \\ 0&\frac{\sqrt{3}}{2} \end{bmatrix} \zed^2$. These observations imply the following conditions on $k$ and $t$:
\begin{equation}
\label{conditions}
\left. \begin{array}{ll}
\sqrt{k} \left( (m+n) \cos t \mp (m-n) \sqrt{3} \sin t \right) \in \zed \\
\sqrt{k} \left( (m-2n) \cos t \mp m \sqrt{3} \sin t \right) \in \zed \\
\sqrt{k} \left( (m+n) \sin t \pm (m-n) \sqrt{3} \cos t \right) \in \sqrt{3} \zed \\
\sqrt{k} \left( (m-2n) \sin t \pm m \sqrt{3} \cos t \right) \in \sqrt{3} \zed
\end{array}
\right\}.
\end{equation}
Let us write $k=3^u j^2 d$, where $j,d \in \zed_{>0}$ with $d$ squarefree and not divisible by 3, and $u=0,1$. We consider two cases.

{\it Case 1.} Suppose first that $u=0$. Then \eqref{conditions} implies that
\begin{equation}
\label{cs_1}
\cos t = \frac{p}{q \sqrt{d}},\ \sin t = \frac{r}{q \sqrt{d}} \sqrt{3}
\end{equation}
for some $q \mid j$.

{\it Case 2.} Suppose next that $u=1$. Then \eqref{conditions} implies that
\begin{equation}
\label{cs_2}
\sin t = \frac{p}{q \sqrt{d}},\ \cos t = \frac{r}{q \sqrt{d}} \sqrt{3}
\end{equation}
for some $q \mid j$.

\noindent
In both cases, the triple $(p,r,q)$ must be a solution to the equation
\begin{equation}
\label{prq3d}
p^2+3r^2=dq^2.
\end{equation}
Notice that \eqref{prq3d} has integral solutions if and only if $d$ is representable by the positive definite binary quadratic form $x^2+3y^2$, in which case all such solutions are given by the parameterization of Lemma \ref{dioph} with $\alpha=1,\beta=0,\gamma=3,\delta=d$. Now, it is a well known fact (see, for instance \cite{butcher}) that $d$ is representable by $x^2+3y^2$ if and only if its prime factorization contains only primes congruent to 1 mod 6 (which, for primes, is the same as $\equiv 1 (\md 3)$); moreover, the number of such representations for a given $d$ is $2^{\omega(d)+1}$, where $\omega(d)$ is the number of distinct prime divisors of $d$, which can be obtained as an easy consequence of unique factorization into irreducibles in the ring of Eisenstein integers. This completes the proof of the lemma.
\endproof

An immediate corollary of this result is an explicit description of the set of all possible index values of WR sublattices of $\Lambda_h$.

\begin{cor} \label{index_values} Let $\II$ be the set of all possible values of $|\Lambda_h: \Gamma|$, where $\Gamma \in \WR(\Lambda_h)$. Then
\begin{eqnarray}
\label{iv_1}
\II = \Big\{ 3^u j^2 d (2m-n) n & : & u=0 \text{ or } 1,\ j,d,m,n \in \zed_{>0}, \nonumber \\
& & d \text{ is $1$ or a product of distinct primes} \nonumber \\
& & \equiv 1 (\md 3),\gcd(m,n)=1, 3 \nmid (m+n), 1 \leq \frac{m}{n} \leq 2 \Big\}.
\end{eqnarray}
\end{cor}

\proof
This follows immediately by combining Lemmas \ref{minimal_gamma} and \ref{sim_class_desc}.
\endproof

Theorem \ref{sim_par} now follows by combining Lemmas \ref{norm_param}, \ref{minimal_gamma}, and \ref{gstm} with Corollary~\ref{norm_angle_cor} (also notice references to Lemma \ref{sim_class_desc} and Corollary \ref{index_values} in the statement of Theorem \ref{sim_par}).
\bigskip

\section{Number, minima, and interference of WR sublattices of $\Lambda_h$}
\label{optimize}

In this section we investigate three related questions, which are the analogues of Questions~1,~2,~and~3 of \cite{sloane} for well-rounded sublattices of $\Lambda_h$. Let $\II$ be as in \eqref{iv_1} and let $J \in \II$.

\begin{quest} \label{Q1} Up to similarity, how many WR sublattices of $\Lambda_h$ of index $J$ are there?
\end{quest}

\begin{quest} \label{Q2} Which WR sublattice of $\Lambda_h$ of index $J$ has the greatest minimum?
\end{quest}

\begin{quest} \label{Q3} Which WR sublattice of $\Lambda_h$ of index $J$ has the highest SNR? 
\end{quest}

\noindent
In what follows, we use our parameterization in Theorem \ref{sim_par} to develop algorithmic procedures and obtain experimental data for Questions \ref{Q1}, \ref{Q2}, \ref{Q3}. We will write $\Gamma_{\theta}(m,n)$ for the lattice as in \eqref{minimal_sub} with the specified choices of $m$ and $n$; we also define
\begin{equation}
\label{A_matrices}
A_1(p,r,q,d) = \begin{bmatrix} \frac{p}{q \sqrt{d}} & -\frac{r  \sqrt{3}}{q \sqrt{d}} \\ \frac{r  \sqrt{3}}{q \sqrt{d}} & \frac{p}{q \sqrt{d}} \end{bmatrix},\ A_2(p,r,q,d) = \begin{bmatrix} \frac{p}{q \sqrt{d}} & \frac{r  \sqrt{3}}{q \sqrt{d}} \\ \frac{r  \sqrt{3}}{q \sqrt{d}} & -\frac{p}{q \sqrt{d}} \end{bmatrix},
\end{equation}
and
\begin{equation}
\label{B_matrices}
B_1(p,r,q,d) = \begin{bmatrix} \frac{r  \sqrt{3}}{q \sqrt{d}} & -\frac{p}{q \sqrt{d}} \\ \frac{p}{q \sqrt{d}} & \frac{r  \sqrt{3}}{q \sqrt{d}} \end{bmatrix},\ B_2(p,r,q,d) = \begin{bmatrix} \frac{r  \sqrt{3}}{q \sqrt{d}} & \frac{p}{q \sqrt{d}} \\ \frac{p}{q \sqrt{d}} & -\frac{r  \sqrt{3}}{q \sqrt{d}} \end{bmatrix}.
\end{equation}
\bigskip

\noindent
{\it Considering Question \ref{Q1}:} Let $\N(J)$ be the number of WR sublattices of $\Lambda_h$ of index $J$. Then  $\N(J)$ is equal to the number of distinct representations of $J$ in the form
\begin{eqnarray}
\label{JJ_values}
J & = & 3^u j^2 d (2m-n) n, \text{ where } u=0 \text{ or } 1,\ j,d,m,n \in \zed_{>0}, \nonumber \\ 
& & d \text{ is 1 or a product of distinct primes} \equiv 1 (\md 3), \nonumber \\
& & \gcd(m,n)=1, 3 \nmid (m+n), 1 \leq \frac{m}{n} \leq 2.
\end{eqnarray}
For a given $J$ it is not difficult to count all such representations. For example, if $J=84$ we have
\begin{eqnarray*}
84 & = & 3^0 \times 2^2 \times 1 \times (2 \times 5 - 3) \times 3 \\
& = & 3^1 \times 2^2 \times 7 \times (2 \times 1 -1) \times 1,
\end{eqnarray*}
and so $\N(84)=2$, where the particular sublattices of $\Lambda_h$ of index 84, up to similarity, are
$$2 \Gamma_{\theta}(5,3),\ 2\sqrt{21} B_1(2,1,1,7)\Lambda_h,$$
corresponding to these representations of 84 respectively. As another example, consider $J=1925$:
\begin{eqnarray*}
1925 & = & 3^0 \times 5^2 \times 1 \times (2 \times 9 - 7) \times 7 \\
& = & 3^0 \times 1^2 \times 7 \times (2 \times 18 - 11) \times 11 \\
& = & 3^0 \times 1^2 \times 35 \times (2 \times 8 - 5) \times 5 \\
& = & 3^0 \times 1^2 \times 55 \times (2 \times 6 - 5) \times 5 \\
& = & 3^0 \times 5^2 \times 77 \times (2 \times 1 -1) \times 1,
\end{eqnarray*}
and so $\N(1925)=5$. The particular sublattices of $\Lambda_h$ of index 1925, up to similarity, corresponding to the first two of these representations are
$$5 \Gamma_{\theta}(9,7),\ \sqrt{7} A_1(2,1,1,7) \Gamma_{\theta}(18,11),$$
while the other three are sublattices similar to $\Gamma_{\theta}(8,5)$, $\Gamma_{\theta}(6,5)$, and $5 \Lambda_h$, respectively.
\bigskip

\noindent
{\it Considering Question \ref{Q2}:} Let $\Gamma \in \WR(\Lambda_h)$ be a lattice with $|\Lambda_h : \Gamma| = J$. Let  $\theta \in C_h$ be such that $\Gamma \in \cht$. Then, on the one hand $\det (\Gamma) = J \det (\Lambda_h) = J \sqrt{3}/2$, and on the other $\det (\Gamma) = |\Gamma| \sin \theta$ by the parallelogram rule. Let $m,n$ be coprime positive integers with $n \leq m \leq 2n$ and $3 \nmid (m+n)$ which correspond to the angle $\theta$ under the parameterization of \eqref{cos_th}, then
$$\sin \theta = \frac{\sqrt{3} (2m-n) n}{2 (n^2-nm+m^2)},$$
and so
\begin{equation}
\label{min_index}
|\Gamma| = \frac{J (n^2-nm+m^2)}{(2m-n)n} \in \zed.
\end{equation}
Hence $(2m-n)n$ must be a divisor of $J (n^2-nm+m^2)$, which implies that $(2m-n)n \mid J$ (since by Lemma \ref{norm_param} $\gcd((n^2-nm+m^2), (2m-n)n) = 1$). Therefore our integers $m,n$ must satisfy the following conditions:
\begin{enumerate}
\item $m,n > 0$ and $\gcd(m,n) = 1$
\item $1 \leq \frac{m}{n} \leq 2$
\item $3 \nmid (m + n)$
\item $(2m-n)n \mid J$
\item $\frac{J (n^2-nm+m^2)}{(2m-n)n}$ is representable by the quadratic form $x^2-xy+y^2$,
\end{enumerate}
and, to maximize $|\Gamma|$, among all such $m,n$ we want to choose a pair that maximizes the quotient $(n^2-nm+m^2) / (2m-n)n$. We can write $m=\beta n$, where $\beta \in [1,2]$, and define
$$f(\beta) = \frac{n^2-nm+m^2}{(2m-n)n} = \frac{\beta^2-\beta+1}{2\beta-1},$$
which is a decreasing function on the interval $[1,(1+\sqrt{3})/2)$ and an increasing function on the interval $((1+\sqrt{3})/2, 2]$. In fact, $f(\beta)$ reaches its maximum at the endpoints of the interval, $\beta=1,2$, which are achieved by the pairs $(m,n)=(1,1)$ and $(2,1)$, respectively. The choice $(m,n)=(2,1)$, however, does not satisfy the condition (3) above, while the choice $(m,n)=(1,1)$ corresponds to the similarity class $\ch(\pi/3) = \left< \Lambda_h \right>$. This means that whenever $J$ is representable by the quadratic form $x^2-xy+y^2$, then $|\Gamma|$ reaches its maximum on $\WR(\Lambda_h)$ at an ideal sublattice of index $J$; this conclusion is consistent with Theorem~3 of \cite{sloane}. As indicated in (12) of \cite{sloane}, this happens for all values of $J$ with prime factorization of the form
\begin{equation}
\label{J_values}
J = 3^k \prod_{p_i \equiv 1 (\md 3)} p_i^{l_i} \prod_{q_j \equiv -1 (\md 3)} q_j^{2m_j},
\end{equation}
where $k,l_i,m_j \in \zed_{\geq 0}$. On the other hand, there exist WR sublattices of $\Lambda_h$ of index $J$ for many values of $J$ not in the form \eqref{J_values}; in these situations, verifying conditions (1)-(5) above for every divisor of $J$ presents a finite search algorithm for the similarity class containing a WR sublattice $\Gamma$ of $\Lambda_h$ of index $J$ with maximal $|\Gamma|$. The next lemma allows to eliminate some of the values of the index $J$ for which WR sublattices of $\Lambda_h$ do not exist.

\begin{lem} \label{eliminate} Suppose that $J \in \zed_{>0}$ is not of the form \eqref{J_values}, and satisfies one of the following:
\begin{enumerate}
\item $J$ is a prime
\item $J=pq$, where $p,q$ are odd primes with $q>3p$
\item $J=2p$, where $p$ is an odd prime
\end{enumerate}
Then $\Lambda_h$ does not have a WR sublattice of index $J$.
\end{lem}

\proof
Suppose there exists $\Gamma \in \WR(\Lambda_h)$ with $|\Lambda_h:\Gamma| = J$, where $J$ is not of the form \eqref{J_values}. Then $\Gamma \in \left< \Gamma_{\theta} \right>$ for some $\Gamma_{\theta} \nsim \Lambda_h$ as in Theorem \ref{sim_par}, meaning that
$$J = k(2m-n)n$$
for some $m,n,k \in \zed_{>0}$ with $\gcd(m,n) = 1$, $n \leq m \leq 2n$, $3 \nmid (m + n)$, and $m,n \neq 1$.

First suppose that $J$ is a prime, then $k=1$, $n=J$, and $m=\frac{J+1}{2}$ (since $n \neq 1$). In this case, however, $\frac{J+1}{2} < J$, contradicting the condition $n \leq m$; hence this is impossible.

Next assume that $J=pq$, where $p,q$ are odd primes with $q>3p$. If $k=p$ or $q$, then $(2m-n)n$ is a prime, and so we are back in the situation above. Then $k=1$, $n=p$, and $2m-p=q$, meaning that $m = (p+q)/2 > 2p = 2n$, which contradicts the condition $m \leq 2n$; hence this is impossible.

Finally, suppose that $J=2p$, where $p$ is an odd prime. If $k=2$ or $p$, then $(2m-n)n$ is a prime, and so we are back in the situation above. Then $k=1$, $n=2$, and $2m-2=p$, which is not possible since $p$ is odd.
\endproof

\noindent
Using Lemma \ref{eliminate}, we can immediately eliminate many of the small values of the index, such as $J=2,5,6,10,11,14,17,22,23,26,29,33,34,38,41,46,47,53,59$, etc. In the table below, we exhibit a few examples we computed for some small values of $J$ not in the form \eqref{J_values}, for which WR sublattices of $\Lambda_h$ of index $J$ exist.

\begin{center} 
\begin{table}[h!b!p!]
\caption{Examples of $\Gamma \in \WR(\Lambda_h)$ of fixed index $J$ maximizing $|\Gamma|$.} 
\begin{tabular}{|l|l|l|l|} \hline
{\em $J = |\Lambda_h : \Gamma|$} & {Maximal $|\Gamma|$} & {Lattice $\Gamma$} \\ \hline \hline
8 & 7 & $\Gamma_{\theta}(3,2)$ \\ \hline
15 & 13 & $\Gamma_{\theta}(4,3)$ \\ \hline
21 & 19 & $\Gamma_{\theta}(5,3)$ \\ \hline
24 & 21 & $\sqrt{3} B_1(1,1,2,1) \Gamma_{\theta}(3,2)$ \\ \hline
32 & 28 & $2 \Gamma_{\theta}(3,2)$ \\ \hline
35 & 31 & $\Gamma_{\theta}(6,5)$ \\ \hline
40 & 37 & $\Gamma_{\theta}(7,4)$ \\ \hline
45 & 39 & $\sqrt{3} B_1(1,1,2,1) \Gamma_{\theta}(4,3)$ \\ \hline
55 & 49 & $\Gamma_{\theta}(8,5)$ \\ \hline
60 & 52 & $2 \Gamma_{\theta}(4,3)$ \\ \hline
65 & 61 & $\Gamma_{\theta}(9,5)$ \\ \hline
\end{tabular}
\label{table1}
\end{table}
\end{center}

\noindent
In fact, there exists an easy test to check our work: if $\Gamma \in \WR(\Lambda_h)$ with $|\Lambda_h:\Gamma| = J$ and $|\Gamma| = M$, then $\Gamma \in \cht$ with
\begin{equation}
\label{cos_check}
\cos \theta = \frac{\sqrt{4M^2 -3J^2}}{2M} \in \que.
\end{equation}
It is now easy to verify that for any value $M$ of $|\Gamma|$ larger than those in the table above, the expression in \eqref{cos_check} would not be rational.
\bigskip

\noindent
{\it Considering Question \ref{Q3}:} As above, suppose that $\Gamma \in \cht$ for some $\theta \in C_h$ is a lattice with $|\Lambda_h : \Gamma| = J$ where $J$ is fixed. Recall that the total interference of $\Gamma$ is given by $E_{\Gamma}(2)$, where $E_{\Gamma}(s)$ is as in \eqref{epstein_z}, and $\SNR(\Gamma)$ is defined in \eqref{SNR_def}. Here we show that Question~\ref{Q3} can be reduced to Question~\ref{Q2}. Notice that this is not generally so for {\it any} (not necessarily WR) sublattices of $\Lambda_h$, as indicated in \cite{sloane}.

\begin{lem} \label{SNR_min} A WR sublattice $\Gamma \subseteq \Lambda_h$ maximizes $\SNR(\Gamma)$ among all WR sublattices of $\Lambda_h$ of index $J$ if and only if it maximizes $|\Gamma|$.
\end{lem}

\proof
Let $M=|\Gamma|$, and let $Q_{\Gamma}(x,y)$ be the quadratic form of $\Gamma$ corresponding to a minimal basis, then
$$Q_{\Gamma}(x,y) = M (x^2+y^2+2xy\cos \theta).$$
The Epstein zeta-function of $\Gamma$ is then given by
$$E_{\Gamma}(s) = \sum_{x,y \in \zed \setminus \{0\}} Q_{\Gamma}(x,y)^{-s} = \frac{1}{M^s} g_{\Gamma}(\theta),$$
where $g_{\Gamma}(\theta) = \sum_{x,y \in \zed \setminus \{0\}} \frac{1}{(x^2+y^2+2xy\cos \theta)^s}$. Then
\begin{eqnarray*}
\frac{d}{d (\cos \theta)} g_{\Gamma}(\theta) & = & \sum_{x,y \in \zed \setminus \{0\}} \frac{2sxy}{(x^2+y^2+2xy\cos \theta)^{s+1}} \\
& = & \sum_{x,y \in \zed_{>0}} \left( \frac{4sxy}{(x^2+y^2+2xy\cos \theta)^{s+1}} - \frac{4sxy}{(x^2+y^2-2xy\cos \theta)^{s+1}} \right) \\
& < & 0,
\end{eqnarray*}
meaning that $g_{\Gamma}(\theta)$ is a decreasing function of $\cos \theta$ for $\theta \in [\pi/3,\pi/2)$. Now $\cos \theta$ is given by \eqref{cos_check}, and is easily seen to be an increasing function of $M$. Therefore for all real values of $s > 1$, $E_{\Gamma}(s)$ is a decreasing function of $M$, in particular meaning that the total interference of $\Gamma$ is minimized (and hence $\SNR(\Gamma)$ is maximized) if and only if $|\Gamma|$ is maximized. 
\endproof

Lemma~\ref{SNR_min} now implies that in order to find a WR sublattice of $\Lambda_h$ of fixed index that maximizes SNR, we can follow the algorithmic procedure described in our consideration of Question~\ref{Q2} above. 
\bigskip

\section{Combinatorial structure in the Eisenstein triples}
\label{combin}

As a part of the parameterization of similarity classes of WR sublattices of $\zed^2$ in \cite{me:wr1}, it has been shown that the set of these similarity classes has the structure of a non-commutative monoid generated by an infinite family of matrices from $\GL_3(\zed)$. The corresponding matrices can be characterized as words of a certain shape in three particular matrices that can be used to generate all primitive Pythagorean triples from $(3,4,5)$ by left multiplication. The existence of these generating matrices for primitive Pythagorean triples has long been known (see \cite{me:wr1} for details). In this section we will investigate a similar generating family of matrices for primitive Eisenstein triples, and use them to explore the combinatorial structure of the set of similarity classes of WR sublattices of $\Lambda_h$.

Let $P_E$ be the set of all primitive Eisenstein triples as defined in section~\ref{norm_f} above. Let also  $G_E = \langle U, M_1, M_2, M_3 \rangle$ be the non-commutative monoid generated by the $3 \times 3$ matrices $U, M_1, M_2, M_3 \in \GL_3(\ZZ)$, defined as follows:
$$      U = 
        \begin{bmatrix} 
         -1&1&0 \\ 
          0&1&0 \\ 
          0&0&1 
        \end{bmatrix},\ 
      M_1 =
        \begin{bmatrix}
         3 & -4 & 4 \\ 
         7 & -7 & 8 \\ 
         6 & -6 & 7
        \end{bmatrix},\ 
      M_2 =
        \begin{bmatrix}
           -4 & 3 & 4 \\
           -7 & 7 & 8 \\
           -6 & 6 & 7
        \end{bmatrix},\ 
       M_3 =
        \begin{bmatrix}
          1 & 3 & 4 \\
          0 & 7 & 8 \\
          0 & 6 & 7
        \end{bmatrix}.$$

\begin{lem} \label{monoid} $G_E$ acts on the set $P_E$ by left multiplication:
\begin{equation}
\label{M_mult}
M(a,b,c) = M \begin{bmatrix} a \\ b \\ c \end{bmatrix},
\end{equation}
for every $M \in G_E$ and $(a,b,c) \in P_E$.
\end{lem}

\proof
A direct verification (for instance, using Maple or SAGE mathematical software packages) shows that for every $1 \leq i \leq 3$ and $(a,b,c) \in P_E$, $M_i(a,b,c) \in P_E$; also $U(a,b,c)=(b-a,b,c) \in P_E$. Since $G_E$ is generated by $U,M_1,M_2,M_3$, it follows that for every $M \in G_E$ and $(a,b,c) \in P_E$, we have $M(a,b,c) \in P_E$.
\endproof

\begin{rem} \label{matrix_conj} Computational evidence (using SAGE) suggests that perhaps all primitive Eisenstein triples can be obtained in this way, starting from $(0,1,1)$; in other words, it seems likely that for each $(a,b,c) \in P_E$ there exists $M \in G_E$ such that $M(0,1,1)=(a,b,c)$. One matrix from the monoid $G_E$ has previously been found in \cite{barnes}.
\end{rem}

The action of Lemma \ref{monoid} induces an action on the set of similarity classes of WR sublattices of $\Lambda_h$. The similarity classes however are parameterized not by the set $P_E$ of primitive Eisenstein triples $(a,b,c)$, but by the set of associated pairs of these triples $\left< a,b,c \right>$ as defined in section~\ref{norm_f} (see Lemma \ref{norm_param} and definition right before it); let us write $\PP_E$ for this set. Recall that the elements of an associated pair are related by $(b-a,b,c) = U(a,b,c)$. Hence we can think of the set $\PP_E$ as the set $P_E$ modulo the equivalence relationship: $(a',b',c') \sim (a,b,c)$ when $(a',b',c') = U(a,b,c)$. For each associated pair in $\PP_E$, let us call the triple with $b>2a$ the {\it upper} triple, denoted $(a,b,c)_u$, and the one with $b<2a$ the {\it lower} triple, denoted $(a,b,c)_l$. Now let us write $M_4 = M_1U, M_5 = M_2U$ and define a sub-semigroup of $G_E$,  $G'_E = \langle M_1, M_2, M_3, M_4, M_5 \rangle$. For each $M \in G'_E$, $\left< a,b,c \right> \in \PP_E$ define 
\begin{equation}
\label{lower_def}
M \left< a,b,c \right> = \left< M(a,b,c)_u \right>,
\end{equation}
where $M(a,b,c)_u$ is as in \eqref{M_mult} for the vector corresponding to the upper triple of $\left< a,b,c \right>$, and so $M \left< a,b,c \right>$ is the associated pair of the triple corresponding to the vector $M(a,b,c)_u$ in $\PP_E$.

\begin{lem} \label{monoid1} $G'_E$ acts on $\PP_E$ by $M \left< a,b,c \right>$ for every $M \in G'_E$, $\left< a,b,c \right> \in \PP_E$.
\end{lem}

\proof
This follows immediately by combining Lemma \ref{monoid} with definition \eqref{lower_def}.
\endproof

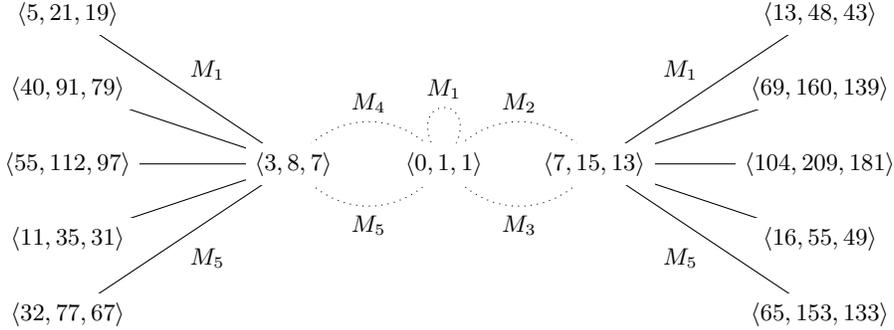
\begin{figure}[h]
        \centering
        \begin{tikzpicture}
             
             % Loop M_1
             \draw (0,0) [dotted] .. controls (-.75,1) and (.75,1) .. node[above] {\small $M_1$} (0,0);
                   
            % (3,8,7) double
            \draw (0,0) [dotted] .. controls (-.5,.75) and (-1.5,.75) .. node[above] {\small $M_4$} (-2,0);
            \draw (0,0) [dotted] .. controls (-.5,-.75) and (-1.5,-.75) .. node[below] {\small $M_5$} (-2,0);
            
            % (7.15,13) double
            \draw (0,0) [dotted] .. controls (.5,.75) and (1.5,.75) .. node[above] {\small $M_2$} (2,0);
            \draw (0,0) [dotted] .. controls (.5,-.75) and (1.5,-.75) .. node[below] {\small $M_3$} (2,0);
            
            % Children
            \draw (-2,0) -- node[above right] {\small $M_1$} (-5,2);
            \draw (2,0) -- node[above left] {\small $M_1$} (5,2);
            \draw (-2,0) -- (-5,1);
            \draw (-2,0) -- (-5,0);
            \draw (-2,0) -- (-5,-1);
            \draw (2,0) -- (5,1);
            \draw (2,0) -- (5,0);
            \draw (2,0) -- (5,-1);
            \draw (-2,0) -- node[below right] {\small $M_5$} (-5,-2);
            \draw (2,0) -- node[below left] {\small $M_5$} (5,-2);
            
            % Nodes on the upper graph
            \draw (0,0) node[fill=white] {\small $\left< 0,1,1 \right>$};
                % Left tree
                \draw (-2,0) node[fill=white] {\small $\left< 3,8,7 \right>$};
                \draw (-5,2) node[fill=white] {\small $\left< 5,21,19 \right>$};
                \draw (-5,1) node[fill=white] {\small $\left< 40,91,79 \right>$};
                \draw (-5,0) node[fill=white] {\small $\left< 55,112,97 \right>$};
                \draw (-5,-1) node[fill=white] {\small $\left< 11,35,31 \right>$};
                \draw (-5,-2) node[fill=white] {\small $\left< 32,77,67 \right>$};
                % Right tree
                \draw (2,0) node[fill=white] {\small $\left< 7,15,13 \right>$};
                \draw (5,2) node[fill=white] {\small $\left< 13,48,43 \right>$};
                \draw (5,1) node[fill=white] {\small $\left< 69,160,139 \right>$};
                \draw (5,0) node[fill=white] {\small $\left< 104,209,181 \right>$};
                \draw (5,-1) node[fill=white] {\small $\left< 16,55,49 \right>$};
                \draw (5,-2) node[fill=white] {\small $\left< 65,153,133 \right>$};
                                    
        \end{tikzpicture}
        \caption{Structure of the set $\PP_E$ induced by the action of the monoid $G'_E$} \label{fig_tree}
    \end{figure}

Computational evidence (using SAGE) suggests that there are no relations between the generators $M_1,\dots,M_5$ of $G'_E$, and that every element of $\PP_E$ can be obtained from $\left< 0,1,1 \right>$ by the action of some $M \in G'_E$. In fact, let $\gamma_1$, $\gamma_2$ be the larger and the smaller roots of the polynomial $143t^2 - 252t + 111$, respectively, and for each $(a,b,c)_u \in P_E$ with $a>0$ define
\[ (x,y,z) = \left\{ \begin{array}{ll}
M_1^{-1}(a,b,c)_u & \mbox{if $\gamma_1 < c/b$} \\
M_4^{-1}(a,b,c)_u & \mbox{if $7/8 \leq c/b < \gamma_1$} \\
M_5^{-1}(a,b,c)_u & \mbox{if $\gamma_2 < c/b < 7/8$} \\
M_2^{-1}(a,b,c)_u & \mbox{if $13/15 \leq c/b < \gamma_2$} \\
M_3^{-1}(a,b,c)_u & \mbox{if $c/b < 13/15$}.
\end{array}
\right. \]
In each of these cases, computational evidence (using SAGE) suggests that $(x,y,z)$ is a primitive upper Eisenstein triple with $z < c$, meaning that $\PP_E$ is the orbit of $\left< 0,1,1 \right>$ under the action of $G'_E$ (and so $P_E$ is the orbit of $(0,1,1)$ under the action of $G_E$). These observations conjecturally imply a nice combinatorial structure on the set $\PP_E$ (and hence on the set of similarity classes of WR sublattices of $\Lambda_h$) of two quinary trees joint at the roots, as illustrated in Figure \ref{fig_tree}.

\bigskip

{\bf Acknowledgment.} We would like to thank the NSF supported Claremont Colleges REU program, under the auspices of which this work was done during the Summer of 2009. We also want to thank Chris Towse for his helpful comments on the subject of this paper.
\bigskip

%\nocite{*}
\bibliographystyle{plain}  % Here the bibliography 
\bibliography{hex_wr}    % is inserted.

\begin{thebibliography}{10}

\bibitem{esm}
A.~H. Banihashemi and A.~K. Khandani.
\newblock On the complexity of decoding lattices using the {K}orkin-{Z}olotarev
  reduced basis.
\newblock {\em IEEE Trans. Inform. Theory}, 44(1):162--171, 1998.

\bibitem{barnes}
F.~Barnes.
\newblock Pythagorean triples, etc.
\newblock \url{http://pythag.net/triples.html}.

\bibitem{sloane}
M.~Bernstein, N.~J.~A. Sloane, and P.~E. Wright.
\newblock On sublattices of the hexagonal lattice.
\newblock {\em Discrete Math.}, 170(1-3):29--39, 1997.

\bibitem{butcher}
J.~C. Butcher.
\newblock Mathematical {M}iniature 9: Hardy's taxi, $x^2 + 3y^2 = p$ and
  {M}ichael {L}ennon.
\newblock {\em New Zealand Math. Soc. Newslett.}, 76:18--19, 1999.

\bibitem{conway}
J.~H. Conway and N.~J.~A. Sloane.
\newblock {\em Sphere Packings, Lattices, and Groups}.
\newblock Springer-Verlag, {T}hird edition, 1999.

\bibitem{me:hex}
L.~Fukshansky.
\newblock Revisiting the hexagonal lattice: on optimal lattice circle packing.
\newblock {\em to appear in Elem. Math.}

\bibitem{me:wr}
L.~Fukshansky.
\newblock On distribution of well-rounded sublattices of {$\zed^2$}.
\newblock {\em J. Number Theory}, 128(8):2359--2393, 2008.

\bibitem{me:wr1}
L.~Fukshansky.
\newblock On similarity classes of well-rounded sublattices of {$\zed^2$}.
\newblock {\em J. Number Theory}, 129(10):2530--2556, 2009.

\bibitem{me:sinai}
L.~Fukshansky and S.~Robins.
\newblock Frobenius problem and the covering radius of a lattice.
\newblock {\em Discrete Comput. Geom.}, 37(3):471--483, 2007.

\bibitem{martinet}
J.~Martinet.
\newblock {\em Perfect Lattices in Euclidean Spaces}.
\newblock Springer-Verlag, 2003.

\bibitem{mcmullen}
C.~McMullen.
\newblock Minkowski's conjecture, well-rounded lattices and topological
  dimension.
\newblock {\em J. Amer. Math. Soc.}, 18(3):711--734, 2005.

\end{thebibliography}
\end{document}